\documentclass{style}

\begin{document}

 \centerline{\bf Coloring $n$-String Tangles}

\centerline{\footnotesize ISABEL K. DARCY\footnote{This work was
supported by a grant from the Joint DMS/NIGMS Initiative to
Support Research in the Area of Mathematical Biology  (NIH
GM76242).} }

\baselineskip=12pt 
\centerline{\footnotesize\it
Department of Mathematics, University of Iowa } \baselineskip=10pt
\centerline{\footnotesize\it  14 MacLean Hall
              Iowa City, IA 52242}

\vspace*{10pt}

\centerline{\footnotesize JUNALYN NAVARRA-MADSEN}
\baselineskip=12pt 
\centerline{\footnotesize\it
Department of Mathematics and Computer Science} \baselineskip=10pt
\centerline{\footnotesize\it Texas Woman's University}
\baselineskip=10pt \centerline{\footnotesize\it MCL 302E, 311
University Drive, Denton, TX 76204}

\begin{abstract}
This expository paper describes how the knot invariant Fox coloring can be applied to tangles. {}{}{}
\end{abstract}


\section{Introduction}

This expository paper describes how the knot invariant Fox coloring
\cite{F1, P1} can be applied to tangles. An $n$-string
tangle $\bf T$ is a 3-ball with $n$ strings
properly
embedded in the 3-ball. The boundary of the 3-ball
and the $2n$ endpoints of the $n$-strings on the
boundary of the 3-ball are not allowed to move.
Tangles were first used by John H. Conway to
tabulate knots [10]. 

Following the presentation in \cite{L}, we will describe coloring via systems of linear
equations so that only an introductory background in linear algebra will be
needed. Fox coloring is related to many
 beautiful areas in topology.
Our interest in this method of coloring links and tangles is to make this paper
 accessible to non-mathematicians as this method is used computationally
  to solve tangle equations arising from protein-DNA interactions \cite{D_mu}.
Also, this approach can make open problems in this area accessible
to undergraduates.  For example, the results of \cite{SWK,  K1} can be
proved using only this linear algebra definition of Fox
coloring combined with a neat trick of Przytycki \cite{P1}.

We will begin with a brief review on coloring knots/links in section \ref{sec_knots}.
In this section we will provide examples, but no proofs.
 For proofs see \cite{L}.  Most of the proofs for knots/links are also  similar to those for tangles given in section \ref{sec_tangles}.
In section 4, we extend the coloring definition to tangles containing a finite number of circles.  In section 5, we give some formulas for determining these invariants for 3-string braids and 2-string rational tangles.  In sections 7 and 8, we discuss embedding tangles in knots.  We make some concluding remarks in section 9.

\section{Coloring knots and links}
\label{sec_knots}

\begin{figure}[htbp]
\vspace*{13pt}
\centerline{\psfig{file=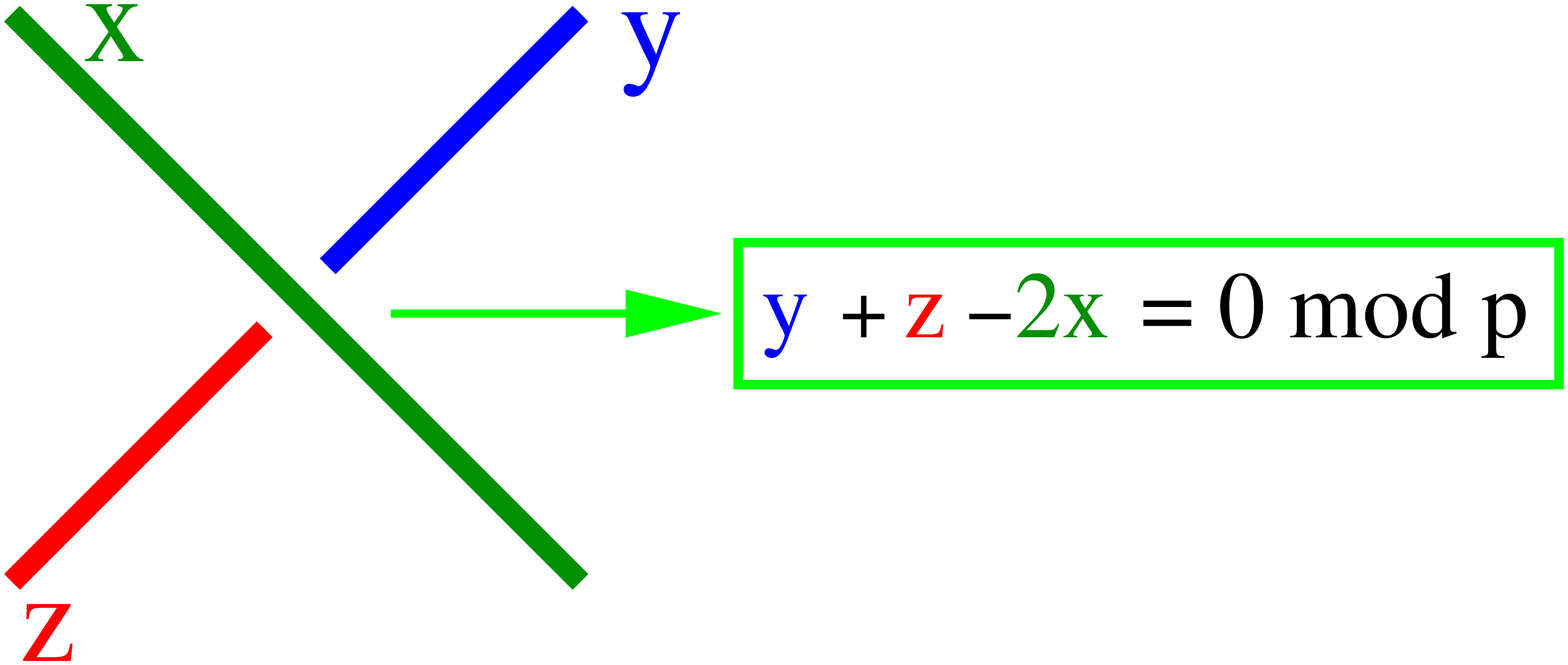,height=.74in,width=1.62in} \quad
\quad \quad \quad \quad \quad \quad \quad
\psfig{file=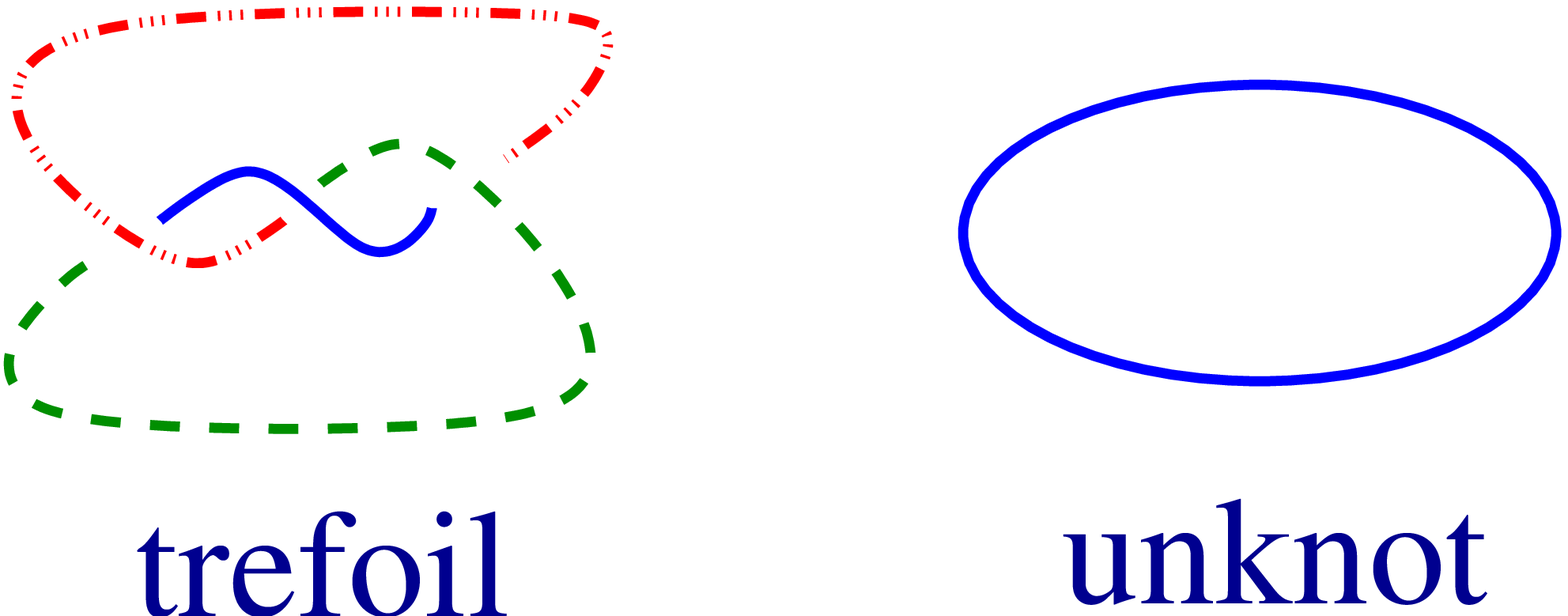,height=.74in,width=1.62in}}
\vspace*{13pt} \caption{A. Coloring condition at a crossing, B.  Trefoil versus unknot} \label{fig:tr_un}
\end{figure}

 An $m$-coloring of a diagram of a knot or
link or tangle is a function $C: \{\it arcs \, of\, a \, diagram\}
\mapsto \mathbb{Z}$$_m$ where the elements of $\mathbb{Z}$$_m =
\{0, 1, \cdots ,$ m-1$\}$ are called {\it colors} and where at at each crossing the 
following relation holds:  if $x$ is the color corresponding to the overarc and $y$ and $z$ are
the two colors corresponding to the two underarcs, then $y+z-2x = 0$  mod $m$ (Fig. 
$\ref{fig:tr_un}$A). 
  If the coloring function is the constant map (i.e., all the arcs are
  assigned the same value or color), then the coloring is said to be {\it trivial}.  
  A link is said to be
  {\it $m$-colorable} if there exists a non-trivial $m$-coloring. 
Coloring mod 3 can easily distinguish a trefoil from an
unknot. Any projection of a trefoil can be colored non-trivially mod 3 while any
projection of an unknot can only be trivially colored. See Fig. $\ref{fig:tr_un}$B.
We explain how to determine if a knot or link is $m$-colorable below using an example.

\vspace{.1cm}

\subsection{Example: A figure-eight knot is 5-colorable}

\begin{figure}[htbp]
\vspace*{13pt}
\centerline{\psfig{file=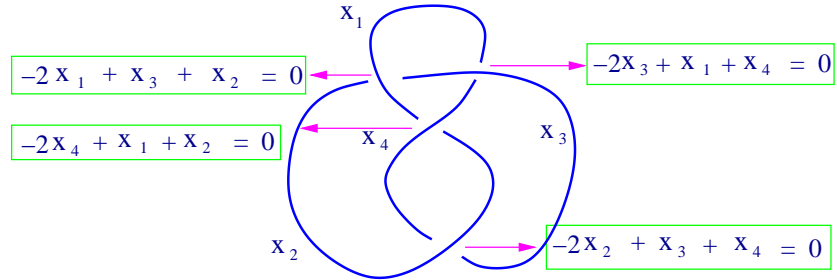,height=1.45in}} 
\vspace*{13pt} \caption{Coloring the figure-eight knot}
\label{fig:4_1}
\end{figure}

Let us color a particular projection of ${\bf 4_1}$ (also called the figure-eight
knot) in figure $\ref{fig:4_1}$. Color each arc of this diagram of ${\bf 4_1}$ using
$x_1, x_2, x_3, x_4$. Note that every knot ${\bf K}$ with $k$
crossings has exactly $k$ arcs. Then, this particular diagram of
${\bf 4_1}$ has four arcs since it has four
crossings in this projection. The crossings are described by the
equations shown in figure $\ref{fig:4_1}$. Writing
these equations in matrix form, we obtain Eqn. \ref{eqn_4_1matrix}

\begin{center}
\begin{equation}
\left(\begin{array}{cccc}
-2 & 1 & 1 & 0 \\
1 & 0 & -2 & 1 \\
0 & -2 & 1 & 1 \\
1 & 1 & 0 & -2 \\
\end{array} \right) \times \left(\begin{array}{c}
x_1\\x_2\\x_3\\x_4\\
\end{array} \right) =  \left(\begin{array}{c}
0\\0\\0\\0\\
\end{array} \right) mod ~m
\label{eqn_4_1matrix}
\end{equation}
\end{center}

Let $\bf M_{4_1} $ be the $4\times 4$ coefficient matrix in Eqn \ref{eqn_4_1matrix}. In order to transform this matrix, $\bf M_{4_1} $,  into echelon form, $EF(\bf M_{4_1})$, we will only use the following elementary row operations:\\

\noindent a) Exchange two rows ($row_i \longleftrightarrow  row_j$)\\
\noindent b) Add a multiple of one row to a different row ($ row_i \longrightarrow
row_i + t \cdot row_j$ where $i \ne j, t \in \mathbb{Z}$)\\
\noindent c) Multiply a row by -1 ($row_i \longrightarrow - row_i)$

\begin{center}
\begin{equation}
EF(\bf{M_{4_1}}) = \left(\begin{array}{cccc}
1 & 0 & -2 & 1 \\
0 & 1 & -3 & 2 \\
0 & 0 & 5 & -5 \\
0 & 0 & 0 &  0 \\
\end{array} \right)
\label{eqn:ef41}
\end{equation}

\end{center}

For those familiar with group presentations, we are forming a finitely generated abelian group where the arcs correspond to generators while the crossing equations give relations among these generators.  Given a knot/link, tangle ${\bf K}$, $\bf{M_{K}}$ is the presentation matrix corresponding to this group.  The allowed row operations allow us to simplify the relations without changing the group.

Each entry of the last row of $EF(\bf M_{4_1})$ (Eqn. \ref{eqn:ef41}) is zero.
The last row of the echelon form of a coloring matrix will always
be a row of all zeros.   This is because all knots and links can
be colored trivially; i.e., given any $a \in {\mathbb{Z}}_m$,
$(x_1,  x_2,\cdots, x_k) = (a, a, \cdots , a)$ will always be a
solution to the system of coloring equations.

From $EF(\bf M_{4_1})$ (Eqn. \ref{eqn:ef41}), we see that if
$m = 5$, the mod 5 solutions to the system of equations in
Eqn. \ref{eqn_4_1matrix} are $(x_1,  x_2, x_3, x_4) = (2a -
b,3a - 2b, a, b)$, $a, b \in {\mathbb{Z}}_5$.
For example if we let $a = 1, b = 0$, then we have the
non-trivial 5-coloring, $(x_1,  x_2, x_3, x_4) = (2,3, 1, 0)$,
Thus ${\bf 4_1}$ is
5-colorable.
If $m$ is a multiple of 5, $m = 5r$ for some $r \in {\mathbb{Z}}$, then
$(x_1,  x_2, x_3, x_4) = (2r,3r, 1r, 0r)$ will be a solution to
$\bf{M_{4_1}}{\bf x} = {\bf 0}$ mod $5r$.  Hence ${\bf 4_1}$ is also
$m$-colorable if $m$ is a multiple of 5 If $m$ is not a multiple of 5,
then 5 is invertible in ${\mathbb{Z}}_m$.  Hence if $m$ is not a
multiple of 5, then ${\bf 4_1}$ can only be trivially colored mod $m$.

Note that we did not use a scaling operation  $row_i \longrightarrow
t \cdot row_j$, $t \in \mathbb{Z}$-$\{ \pm 1 \}$, as part of the three row
operations above to convert ${\bf M_{4_1}}$ to echelon form,
$EF({\bf M_{4_1}})$. In the example above, notice that scaling by $t =
\frac{1}{5}$ was not done on the third row of  $EF({\bf M_{4_1}})$.
Had we scaled by $t = \frac{1}{5}$, we would lose the information
that the linear system in equation $\ref{eqn:ef41}$ has a
nontrivial $ mod\,\,5$ solution. Also, had we scaled the third row
by 3, we would have been led to the false conclusion that ${\bf 4_1}$ is
3-colorable which it is not.

There are other invariants that can be gleaned from this matrix
method of presenting a link. Since all knots/links have $m$ trivial $m$-colorings, the determinant of
the coloring matrix is always zero.  However, the absolute value of the
determinant of the matrix obtained after removing one row and one column,
${\bf d(L)}$, is an invariant.   For example, ${\bf d(4_1)} = 5$.  A link ${\bf L}$ is $m$-colorable if and only if $gcd(m, {\bf d(L)}) > 1$.    The coloring
matrix with one row and one column removed is the same as the Alexander matrix when $t = -1$. Hence this determinant
is actually the Alexander polynomial evaluated at -1.  For more information
on the Alexander matrix/polynomial see \cite{L}.

\section{Coloring of $n$-string Tangles}
\label{sec_tangles}

\begin{figure}[htbp] 
\vspace*{13pt}
\centerline{\psfig{file=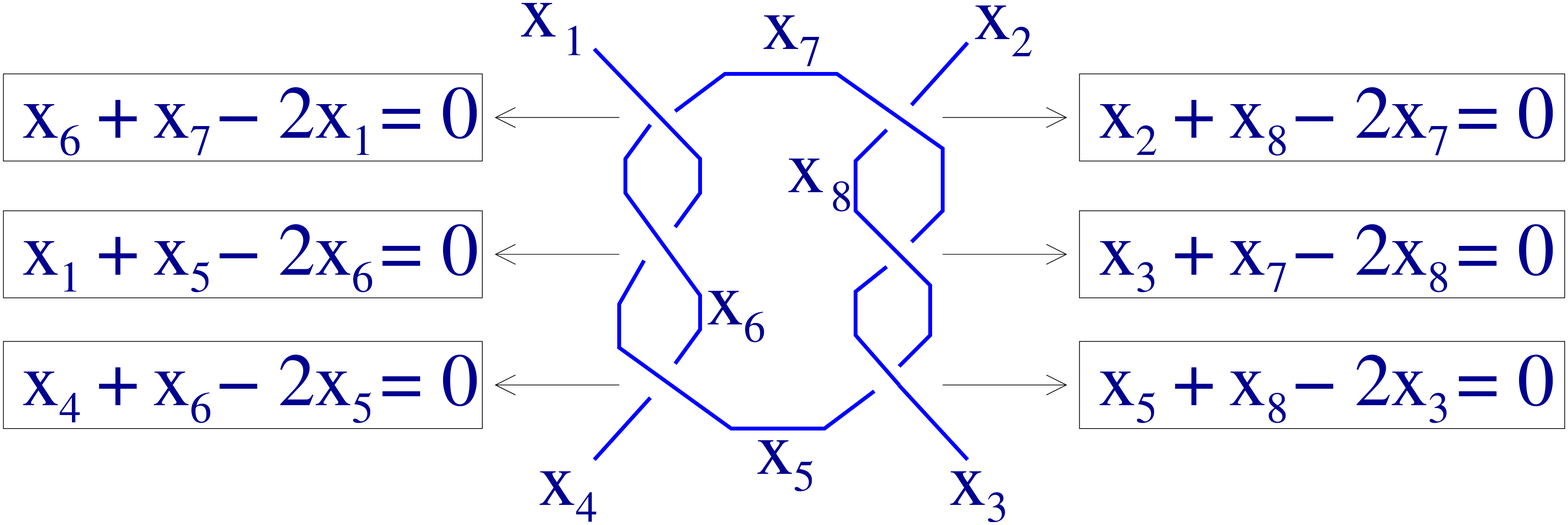,height=1.1in,width=4.1in}} 
{ \vskip 13pt} \caption{Coloring A 2-string Tangle}
\label{fig:bcc}
\end{figure}

We can similarly color $n$-string tangles.  One of the invariants coming from coloring a tangle will depend on  a chosen ordering of the $2n$ endpoints of the $n$ strings. 
We will call the arcs which have one endpoint on the boundary of the 3-ball
{\it endpoint arcs}.  We will fix a particular ordering for the endpoint arcs.
For example, for a 2-string tangle, we will label the endpoint arcs in a clockwise manner
starting with labeling the top left arc $x_1$ as in Fig. \ref{fig:bcc}.
The arcs which are not endpoint arcs will be called {\it interior arcs}.

Note that from Fig. $\ref{fig:bcc}$, at each crossing we form an
equation just as in the knot or link case. This equation
represents a row in the matrix we are going to form out of this
colored tangle. Each arc represents a column of this matrix. A matrix which is row equivalent to a matrix which comes from coloring a tangle ${\bf T}$ will be called  a {\it coloring matrix of ${\bf T}$}.  In the $2$-string tangle example in Fig. $\ref{fig:bcc}$, there are six crossings and eight arcs, thereby
giving a $6 \times (6+2)$ coloring matrix with  entries of zeroes, ones and
negative twos. This is one of the differences between knots and
tangles. We end up having a non-square matrix when coloring a
tangle.  Normally if an $n$-string tangle has $k$ crossings, it will have  $k+n$
arcs, and hence its coloring matrix will have $k$ rows and $k+n$ columns.  See the note just before the proof of Theorem \ref{thm_tangle} for an example of an $n$ string tangle with $k$ crossings for which we choose a coloring matrix which is not $k \times (k+n)$.

Based on the labeling of the given $2$-string tangle in Fig. $\ref{fig:bcc}$, we get the system of linear equations, ${\bf
(M_{T})x} = 0$ in Eqn. \ref{eqn:sq1}.  Notice that we put the endpoint arcs' unknowns, $x_1, x_2, x_3, x_4$,
as the four rightmost columns of matrix $\bf{M_T}$.

\begin{center}
\begin{equation}
\left(\begin{array}{cccccccc}
 0 & 1 & 1 & 0     & -2 & 0 & 0 & 0\\
1 & -2 & 0 & 0     & 1 & 0 & 0 & 0\\
-2 & 1 & 0 & 0    & 0 & 0 & 0 & 1\\
 0 & 0 & -2 & 1    & 0 & 1 & 0 & 0\\
 0 & 0 & 1 & -2    & 0 & 0 & 1 & 0\\
1 & 0 & 0 & 1     & 0 & 0 & -2 & 0\\
\end{array} \right) \times \left(\begin{array}{c}
x_5\\x_6\\x_7\\x_8\\x_1\\x_2\\x_3\\x_4\\
\end{array} \right) =  \left(\begin{array}{c}
0\\0\\0\\0\\0\\0\\0\\0\\
\end{array} \right)
\label{eqn:sq1}
\end{equation}
\end{center}

After
performing the allowed elementary row operations, we obtain an echelon form of
$\bf{M_T}$. Recall that scaling a row is not allowed. An echelon form, $EF(\bf{M_T})$  is:

\begin{center}
\begin{equation}
EF(\bf{M_T}) = \left(
\begin{matrix}
  1 & 0 & 0 & 1     & 0 & 0 & -2 & 0\cr
  0 & 1 & 1 & 0     & -2 & 0 & 0 & 0\cr
  0 & 0 & 1  & -2    & 0 & 0 & 1  & 0\cr
  0 & 0 & 0  & 3    & -3 & 0 & 0  & 0\cr
  0 & 0 & 0  & 0     & 1 & -1 & 1 & -1\cr
  0 & 0 & 0 & 0     & 0 & -2 & 5  & -3\cr
\end{matrix}\right)
\label{eqn_EF1}
\end{equation}
\end{center}

Since we are dealing with $n = 2$ strings, we are interested in the
lower right hand corner $2 \times 4$ submatrix of the echelon form of the matrix $\bf{M_T}$. We show below in general that for an $n$-string tangle, this lower right-hand corner $n
\times 2(n)$ submatrix is an invariant up to the allowed elementary row
operations.  To make it an invariant, we define the standard echelon form of a matrix.

Let $EF(\bf{M_T})=\{a_{ij}\}_{1\leq i \leq k, 1 \leq j \leq
(k+n)}$ be an echelon form of a matrix of size $k \times (k+n)$.  A $leading\,\, entry$ of $EF(\bf{M_T})$ is the first nonzero
entry of a row of $EF(\bf{M_T})$. A matrix $\bf{M_T}$
is in {\it standard echelon form} if the following three properties
hold:

\begin{romanlist}[(ii)]
\item It is in echelon form.
\item Its leading entries are positive.
\item If $a_{ij}$ is a leading entry of the $i$th row, then $0\le a_{\lambda
j} \le a_{ij}-1, 1\le \lambda < i$, i.e., all entries above a
leading entry are non-negative and less than that leading entry.
\end{romanlist}

\begin{lemma}
Let $SF(\bf{M_T})$ be the standard echelon
form of a matrix ${\bf M_T}$. Then $SF(\bf{M_T})$ is unique.
\end{lemma}

\begin{proof}
Similar to showing reduced echelon form is unique.
\end{proof}

From Eqn. \ref{eqn_EF1}, we obtain the standard echelon
form, $SF(\bf{M_T})$, based on the three criteria listed above:
\begin{center}
\begin{equation}
SF(\bf{M_T})= \left(
\begin{matrix}
  1 & 0 & 0 & 1     & 0 & 0 & -2 & 0\cr
  0 & 1 & 1 & 0     & -2 & 0 & 0 & 0\cr
  0 & 0 & 1  & -2    & 0 & 0 & 1  & 0\cr
  0 & 0 & 0  & 3    & -3 & 0 & 0  & 0\cr
  0 & 0 & 0  & 0     & 1 & 1 & -4 & 2\cr
  0 & 0 & 0 & 0      & 0 & 2 & -5 & 3\cr
\end{matrix}\right)
\label{eqn_SF}
\end{equation}
\end{center}

\begin{theorem}
\label{thm_tangle}
Suppose we have chosen a fixed ordering of the endpoint arcs.  The following are invariants of an $n$-string
tangle $\bf{T}$:
\begin{enumerate}
\item ${\bf d_U(T)}$ = absolute value of the the determinant of the upper left
hand corner $(k-n)\times (k-n)$ submatrix  of $\bf T$. 
\item
${\bf
M_{l}(T)}$ = the $n \times 2n$ lower right hand corner submatrix of $SF(\bf{M_T})$. 
\end{enumerate}
\end{theorem}

For the tangle in Fig. \ref{eqn:sq1}, ${\bf d_U(T)}$ = 3 and  ${\bf
M_{l}(T)} = \left(\begin{array}{cccc}
1  & 1  & -4 & 2 \\
0  & 2  & -5 & 3 \\\end{array} \right) $. 
\vskip 10pt

\begin{figure}[htbp] 
\centerline{\psfig{file=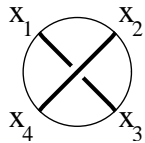,height=.7in}} 
\caption{Arcs can be doubly labeled.}
\label{fig_1tangle}
\end{figure}

\hskip -15pt{\bf{Note:}}
It is possible that one or more strings of an $n$-string tangle does not pass under any string.  
Such a string will project to a single arc.  See, for example, the tangle in Fig. \ref{fig_1tangle}.  
To calculate the invariant ${\bf M_{l}(T)}$ for an $n$-string tangle, {\bf T}, we need to have $2n$ distinct 
variables corresponding to endpoint arcs.  Hence if any string projects to a single arc, we will doubly 
label this arc  with two variables, $x_i$ and $x_j$ (depending on the ordering of the endpoint arcs) 
and add the equation, $x_i - x_j = 0$.  For example, the matrix in Eqn. \ref{eqn_one} is a coloring 
matrix for the one crossing tangle in Fig. \ref{fig_1tangle}.   Additionally, doubly labeling any arc and adding 
an equation(s) equating the variables corresponding to this doubly  labeled arc does not affect 
the invariants listed in Theorem \ref{thm_tangle}.  

\begin{equation}
\left(\begin{array}{cccc}
1  & 0  & 1 & -2 \\
0  & 1  & 0 & -1 \\
\end{array} \right) 
\label{eqn_one}
\end{equation}

\begin{proof}[Theorem \ref{thm_tangle}]

Two coloring matrices of a tangle diagram may differ with respect to how the interior arcs are labeled.  One can convert between $k \times (k + n)$ coloring matrices for the same tangle diagram which differ with respect to interior arc labeling by performing column operations on the first $k - 2n$ columns.  Such column operations only affect the sign of the determinant and do not affect the lower right $n \times 2n$ matrix since no column operations are performed on the last 2$n$ columns.  Similarly doubly labeling an arc has no affect on ${\bf d_U(T)}$ and ${\bf M_{l}(T)}$.

Our allowed row operations can only change the sign of the determinant.    Also, no matter how the allowed row operations are performed, $SF(\bf{M_T})$   is unique.

Hence ${\bf d_U(T)}$ and ${\bf M_{l}(T)}$ are invariants of a given tangle diagram.  Thus we only need to check if they are the same for two different tangle diagrams corresponding to the same tangle.   Hence we only need to check if they are preserved under Reidemeister moves (Fig. \ref{fig_reid}).

\begin{figure}[htbp] 
\centerline{\psfig{file=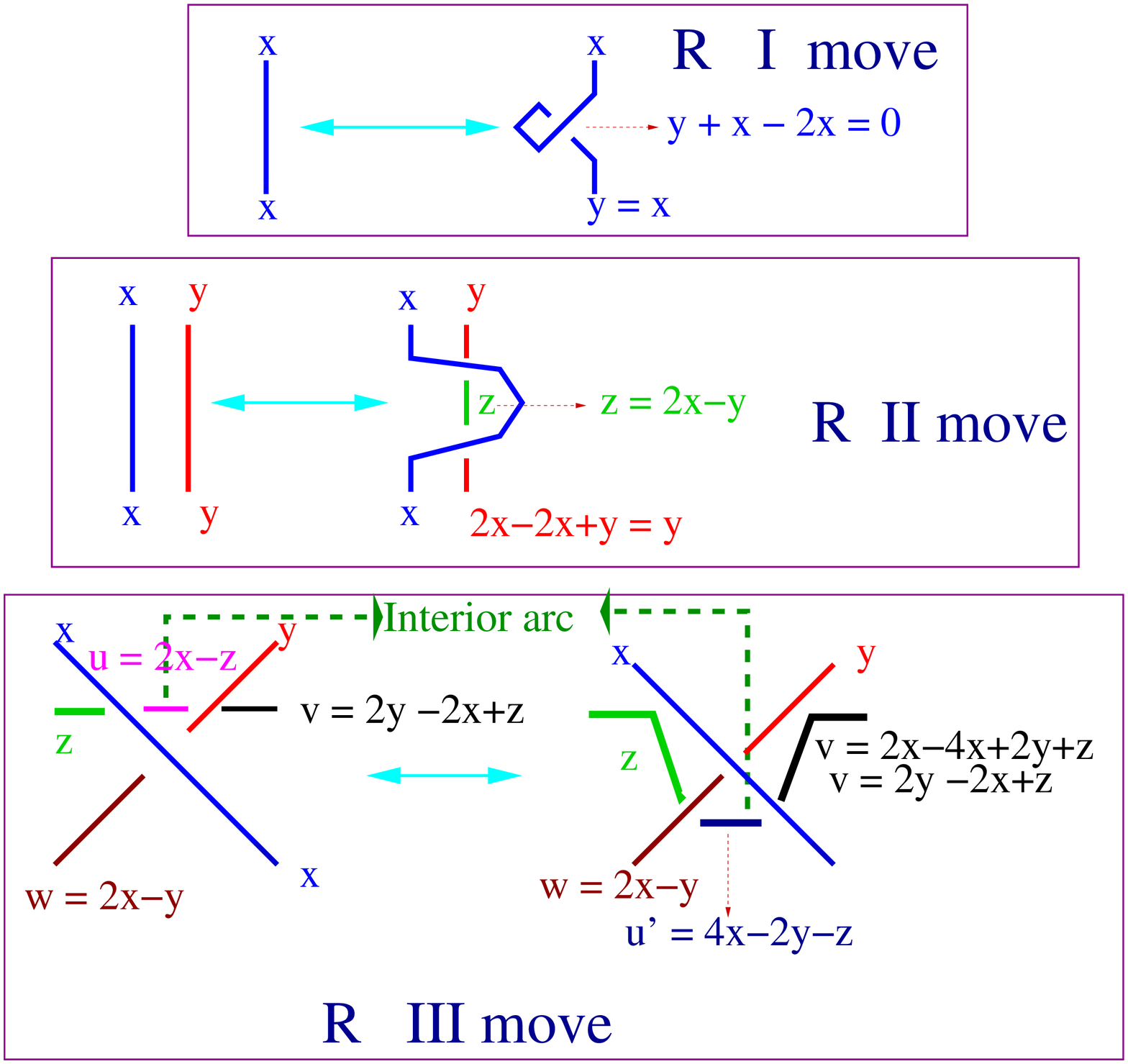,height=3.7in,width=4.4in}} 
\caption{The colors of the endpoint arcs Are Preserved Up To Reidemeister Moves.}
\label{fig_reid}
\end{figure}

A Reidemeister move can be thought of in terms of modifying a subtangle within a tangle.  For example, an RI move can be thought of as replacing a 1-string subtangle containing no crossings with a 1-string subtangle containing exactly one crossing (or vice versa).  
An RI move results in the addition (or removal) of a crossing and the creation (or deletion via joining) of a new arc.  This results in the addition (or removal) of one equation, $x - y = 0$, and one variable. 
The new
equation can be used to eliminate the new variable from all other equations. 
Since $x = y$, the RI move does not affect the color of the endpoint arc(s) of the RI subtangle (Fig. \ref{fig_reid}, top).  
Hence after eliminating the new variable from equations resulting from crossings outside of the RI subtangle, the only difference between the coloring matrices is the addition (or removal) of a row and column. Hence since the endpoint colors are not affected, ${\bf M_l(T)}$ is unchanged by an RI move.
As the leading entry of the added (or deleted) row is 1, the determinant, ${\bf d_U(T)}$, is unchanged by an RI move.

Similarly an RII move consists of modifying a 2-string subtangle 
(Fig. \ref{fig_reid}, middle). In this case, an RII move results in the creation (or deletion) of two new crossings and two new arcs.  Since the endpoint arc colors of the RII subtangle are not affected by an RII move, we can again remove the new subtangle endpoint arc variable from any equation resulting from crossings outside of the RII subtangle so that these equations are identical both before and after the RII move.  Hence 
${\bf M_l(T)}$ is unchanged by an RII move.  Also since the leading entries of the added (or deleted) rows are 1, the determinant, ${\bf d_U(T)}$,  is unchanged by an RII move.

Similarly, the endpoint arc colors of the 3-string RIII subtangle are not affected by an RIII move (Fig. \ref{fig_reid}, bottom).  Hence ${\bf M_l(T)}$ is unchanged by an RIII move.  The equation corresponding to the interior arc of the RIII subtangle is affected by an RIII move, but as an interior arc of the RIII subtangle, this does not affect ${\bf M_l(T)}$.  Since the leading entry of the row corresponding to this interior arc of the RIII subtangle is 1, the determinant, ${\bf d_U(T)}$, is also unchanged by an RIII move.

Thus ${\bf M_{l}(T)}$ and ${\bf d_{U}(T)}$ are not affected by any of the Reidemeister
moves and hence are tangle invariants.
\end{proof}

\section{Other definitions of tangle coloring}

\begin{figure}[htbp] 
\centerline{\psfig{file=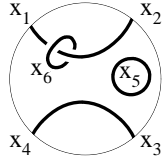,height=0.8in} }  
\caption{A tangle containing two arcs and two circular components} 
\label{fig_tangcircle}
\end{figure}

We defined an $n$-string tangle as  3-ball containing  $n$ properly embedded arcs.  Sometimes one would also like to allow a finite number of circles to be embedded within the 3-ball. We can also apply coloring to these tangles.  In this case we not only label all arcs in the tangle diagram (including those from both strings and circles), but we also label any circular component in the tangle diagram.  We also add a row of all zeros for each such closed circular component in the tangle diagram. For example, the 2-string tangle in Fig. \ref{fig_tangcircle} contains two circles.  One of these projects to a closed  circular component (labeled $x_5$) while the other projects to single arc (labeled $x_6$) in this tangle diagram.  The former results in Eqn. 3: 0 = 0 while the later is involved in two equations,  Eqn. 1:  $x_1 + x_2 - 2x_6 = 0$ and Eqn. 2: $x_6 + x_6 - 2x_2 = 2x_6 - 2x_2 = 0$.  We also have a fourth equation equating two endpoint arcs, $ x_3 - x_4 = 0$.  Hence we obtain the coloring matrix in Eqn. \ref{eqn_circs}A.

  In order to obtain the most information from the coloring equations, we will not use $SF(\bf{M_T})$  in this case.   We can obtain an echelon form, but we will then place  the $n$ rows with a leading entry corresponding to an endpoint arc as the last $n$ rows even below rows of all zero's.  Hence after obtaining an echelon form, all rows of all zeros should be moved above the last $n$ rows.  Thus we obtain the matrix in Eqn. \ref{eqn_circs}B.  Thus $d_u(T) = 0$ and ${\bf M_{l}(T)}= \left(\begin{array}{ cccc }
1 & -1 & 0 & 0 \cr
0 & 0  & 1 & -1 \cr
\end{array} \right)$.

\begin{equation}
{\bf A.)~} \left(
\begin{array}{cccccc }0 & -2 & 1 & 1 & 0 & 0 \cr
0 & 2 & 0 & -2 & 0 & 0 \cr
0 & 0 & 0 & 0  & 0 & 0 \cr
0 & 0 & 0 & 0  & 1 & -1 \cr
\end{array} \right)
\hskip 0.5in
{\bf B.)~}  \left(
\begin{array}{cc|cccc }
0 & 2 & 0 & -2 & 0 & 0 \cr
0 & 0 & 0 & 0  & 0 & 0 \cr \hline
0 & 0 & 1 & -1 & 0 & 0 \cr
0 & 0 & 0 & 0  & 1 & -1 \cr
\end{array} \right)
\label{eqn_circs}
\end{equation}

Jozef Pryztycki determined a relationship among endpoint arcs which all tangles must satisfy:


\vskip 10pt
\begin{theorem}\cite{P1}
If ${\bf T}$ is an $n$-string tangle (possibly containing a finite number of circles), then ${\bf M_l(T)}$ is row equivalent to a matrix where the first row consists of alternating 1's and -1's.
\label{thm_P}
\end{theorem}

\begin{figure}[htbp] 
\centerline{\psfig{file=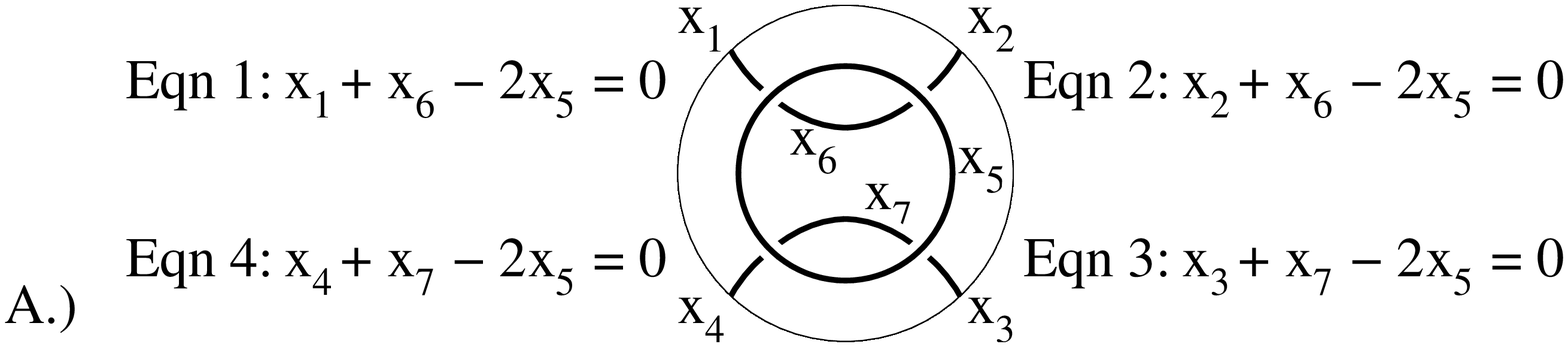,height=0.8in} \hskip 0.1in \hfill \psfig{file=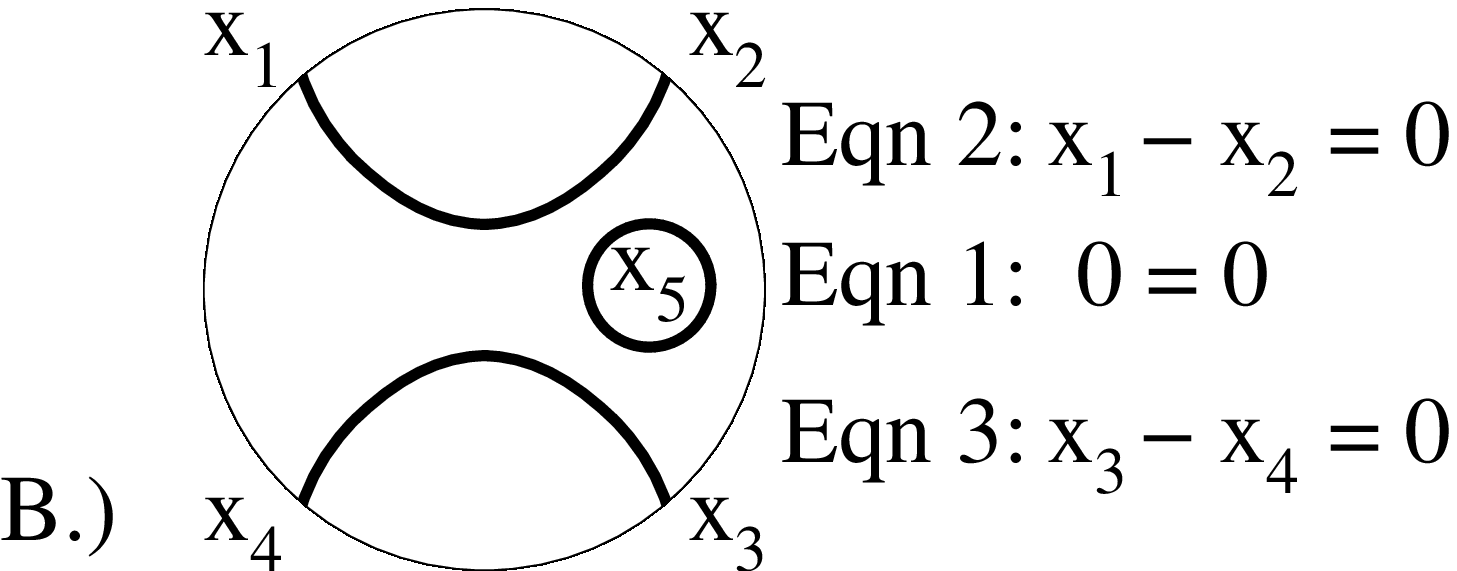,height=0.8in}}  
\caption{} 
\label{fig_circle}
\end{figure}

We will illustrate his theorem and proof with an example.  The tangle in Fig. \ref {fig_circle}A contains an unknotted, unlinked circle.
Note for the tangle diagram in Fig. \ref {fig_circle}A, we have the following relationship among the endpoint arcs:
Eqn 1 - Eqn 2 + Eqn 3 - Eqn 4 = 0.  Hence $(x_1 + x_6 - 2x_5) - (x_2 + x_6 - 2x_5) 
+ (x_3 + x_7 - 2x_5) - (x_4 + x_7 - 2x_5) = x_1 - x_2 + x_3 - x_4 = 0$.
Since both the tangles in Fig. \ref{fig_circle} are the same, the coloring equations corresponding to the tangle diagram in Fig. \ref{fig_circle}B also satisfies the endpoint arcs relationship, $x_1 - x_2 + x_3 - x_4 = 0$.  Removing the circle $x_5$ from the tangle diagram in Fig. \ref{fig_circle}B corresponds to removing the column corresponding to $x_5$ (containing all zeros) as well as the row containing all zeros.  The remaining equations are unchanged.  Hence, ${\bf M_{l}(T)}$ is not affected and we still have the relationship $x_1 - x_2 + x_3 - x_4 = 0$ for the tangle diagram without the circle $x_5$. Note that we can add an unknotted, unlinked circular component to any $n$-string tangle to determine the endpoint arcs relationship, $x_1 - x_2 + ... + x_{n-1} - x_n = 0$.  Since removing the unknotted, unlinked component does not affect ${\bf M_l(T)}$, we can see that the coloring equations of any $n$-string tangles  must satisfy this relationship \cite{P1} thus giving us Theorem \ref{thm_P}.

\begin{equation}
{\bf A.)}\left(
\begin{matrix}
-2 & 1 & 0 & 1 &  0 & 0 & 0  \cr 
-2 & 1 & 0 & 0 &  1 & 0 & 0  \cr 
-2 & 0 & 1 & 0 &  0 & 1 & 0  \cr 
-2 & 0 & 1 & 0 &  0 & 0 & 1  \cr 
 0 & 0 & 0 & 0 &  0 & 0 & 0  \cr
 \end{matrix} \right)
~~\sim~~
\left(\begin{array}{ccc|cccc}
-2 & 0 &  1 & 0 &  0 &  1 & 0  \cr 
 0 & 1 & -1 & 1 &  0 & -1 & 0  \cr 
0 & 0 & 0 & 0 &  0 & 0 & 0  \cr
\hline
 0 & ~0 & 0 & 1 &  -1 & 1 & -1  \cr 
 0 & ~0 & 0 & ~0 &  0 & 1 & -1  \cr 
\end{array} \right)
\hskip 1cm  \hfill
{\bf B.)}
\left(
\begin{array}{c|cccc}
0 & 0 & 0  & 0 & 0 \cr \hline
0 & 1 & -1 & 0 & 0 \cr
0 & 0 & 0  & 1 & -1 \cr
\end{array} \right)
~~\sim~~
\left(
\begin{array}{c|cccc}
0 & 0 & 0  & 0 & 0 \cr \hline
0 & 1 & -1 & 1 & -1 \cr
0 & 0 & 0  & 1 & -1 \cr
\end{array} \right)
\end{equation}

Hence instead of using $SF(\bf{M_T})$, we will sometimes use a matrix row equivalent to ${\bf M_l(T)}$ where the first row consists of alternating 1's and -1's.

\section{3-string braid  and 2-string rational tangle coloring formulas}

\begin{figure}[htbp] 
\centerline{\psfig{file=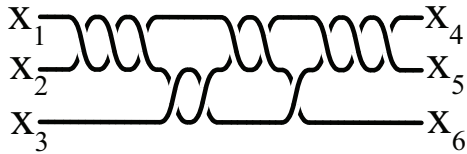,width=2.4in}}
\caption{A 3-string braid}
\label{fig_braid}
\end{figure}

Let us start with a formal definition of an $n$-string braid.  An  $n$-string braid is the union ${\bf B = b_1 \cup b_2 \cup \cdots \cup b_n}$ of $n$ strings $b_i (i = 1, 2, \cdots, n)$ in the cylinder $D^2 \times [0,1]$ such that for each $t \in [0,1]$, ${\bf B}$ intersects the 2-disk $D^2 \times \{t\}$ transversely in $n$ distinct interior points of $D^2 \times \{t\}$ with the $2n$ endpoints fixed. An example of a 3-string braid is given in Fig. \ref{fig_braid}.

In order to calculate coloring formulas for 3-string braids, we will use the Euler bracket function, $E[c_1,..., c_h]$ which equals the sum of products of the $x_i$'s where zero or more disjoint
pairs of consecutive $x_i$'s are omitted \cite{Roberts}. 
For example, $E[c_1, c_2] = c_1c_2+ 1$, $E[c_1, c_2, c_3] = c_1c_2c_3 + c_1 + c_3$, 
$E[c_1, c_2, c_3, c_4] = c_1c_2c_3c_4 + c_1c_2 + c_1c_4 + c_3c_4 + 1$.
If $h = 0$,
then $ E[] = 1$.  Two useful formulas involving the Euler bracket are
$E[c_1,..., c_h] = E[c_h,..., c_1]$
and
$E[c_1,..., c_h] = c_1 E[c_2,..., c_h] + E[c_3,..., c_h]$ \cite{Roberts}.  The following theorem is an unpublished result of Arun Ponnusamy and D.

\begin{theorem}
If ${\bf B}$ is an $n$-string braid, then ${\bf d_U(B)} = 1$.  Furthermore, for a 3-string braid, ${\bf B} = \sigma_1^{-c_1}\sigma_2^{c_2}\sigma_1^{-c_3}...\sigma_2^{c_{h-1}}\sigma_1^{-c_h}$, $h$ odd, if the endpoint arcs have been ordered as in Fig. \ref{fig_braid}, then $M_l({\bf B})$ is row equivalent to the matrix in Eqn. \ref{eqn_braidE}.

\begin{equation}
\left(
\begin{matrix}
   1 & -1 & ~1 &  -1 & 1 & -1 \cr
   0 & ~1 & ~0 &  ~~~ E[1, c_1, ..., c_h] - 1
         ~~~~   &  -E[1, c_1, ..., c_h, 1]    + 1      ~~~~
		   &   E[1, c_1, ..., c_{h-1}] - 1 
\cr
   0 & ~0 & ~1 &   E[c_2,..., c_h] - 1
		   &  -E[c_2,..., c_h, 1] + 1 
		   &   E[c_2..., c_{h-1}] - 1
\end{matrix} \right)
\label{eqn_braidE}
\end{equation}

\end{theorem}

\begin{proof}
Induction on $h$.
\end{proof}

\begin{figure}[htbp]
\centerline{\psfig{file=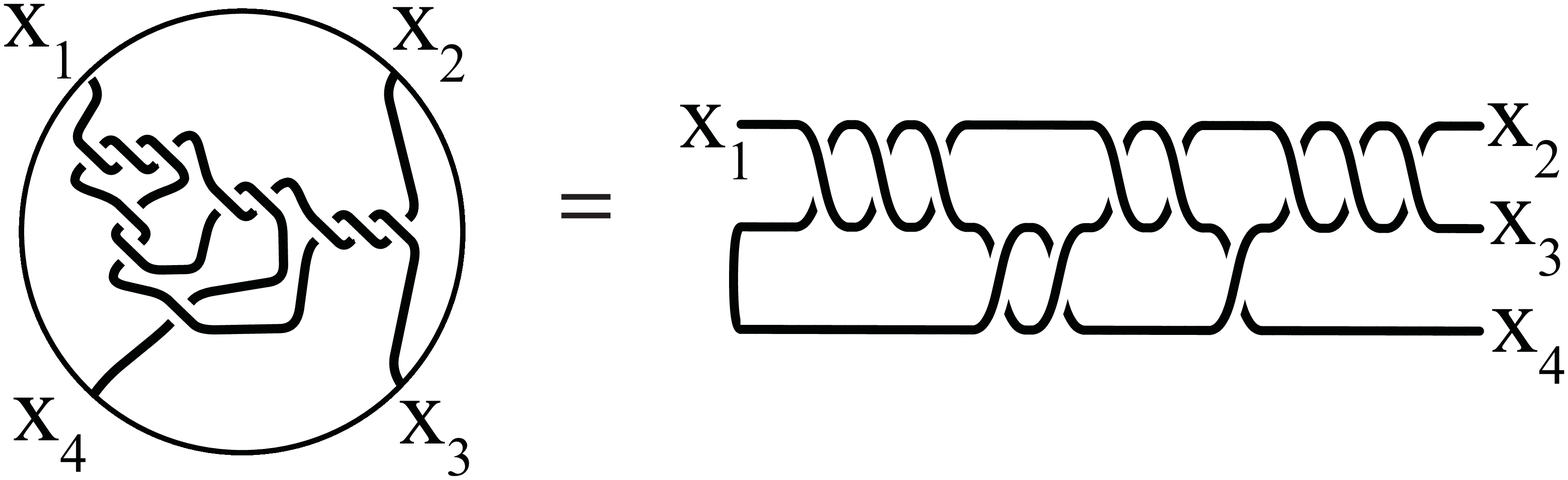,width=3.4in}}
\caption{}
\label{fig_rat}
\end{figure}

A rational 2-string tangle is a tangle which can be formed from a 3-string braid by connecting the endpoint arcs $x_2$ and $x_3$ (compare Figs.
\ref{fig_braid}, \ref{fig_rat}).  For other definitions of rational tangle, see for example \cite{}.  The 3-string braid ${\bf B} = \sigma_1^{-c_1}\sigma_2^{c_2}\sigma_1^{-c_3}...\sigma_2^{c_{h-1}}\sigma_1^{-c_h}$, $h$ odd,
forms the 2-string tangle $<c_1, ..., c_h>$.  2-string tangles are uniquely identified by the continued fraction \cite{Conway, KG1}:  ${p \over q} = {E[c_1, ..., c_h] \over E[c_1, ..., c_{h-1}]}
= c_h + {1 \over {c_{h-1} + ... + {1 \over c_1}}}$.

Thus a coloring matrix for ${\bf T}$ can be obtained from a coloring matrix for ${\bf B}$ by adding the equation $x_2 = x_3$ and switching the columns corresponding to old endpoint arcs of ${\bf B}$, $x_1$ and $x_3$, so that the 2-string tangle endpoint arcs are the last four columns of the coloring matrix (see Eqn. \ref{eqn_rat}).
Hence if ${\bf T}$ is a rational 2-string tangle, then ${\bf d_U(T)} = 1$.  We can determine ${\bf M_l(T)}$ by putting the matrix in Eqn. \ref{eqn_rat} into echelon form.  Hence we have Theorem \ref{thm_rat}.

\begin{equation}
\left(
\begin{matrix}
   1 & -1 & ~1 &  -1 & 1 & -1 \cr
   0 & ~1 & ~0 &  ~~~ E[1, c_1, ..., c_h] - 1
         ~~~~   &  -E[1, c_1, ..., c_h, 1]    + 1      ~~~~
		   &   E[1, c_1, ..., c_{h-1}] - 1 
\cr
   1 & ~0 & ~0 &   E[c_2,..., c_h] - 1
		   &  -E[c_2,..., c_h, 1] + 1 
		   &   E[c_2..., c_{h-1}] - 1
\cr
   1 & -1 &    0 & ~0 &   0 & ~0 
\end{matrix} \right)
\label{eqn_rat}
\end{equation}
\vskip 9pt

\begin{theorem}
\label{thm_rat}
For a rational 2-string tangle ${\bf T} = {\bf {p \over q}}$, ${\bf d_U(T)} = 1$ and 
${\bf M_l(T)} =\left(
\begin{matrix}
1 ~~& -1 ~~& 1      ~~& -1   \cr
0 ~~& p  ~~& -p - q ~~& q
\end{matrix} \right)$
\end{theorem}

Thus, as also noted by \cite{KL}, coloring classifies rational tangles.

\section{Numerator and Denominator Closure of 2-string tangles}

There are a number of operations which can be performed on tangles in order to obtain knots or links.  In this section we will look at the operations of numerator and denominator closures of 2-string tangles.

Suppose ${\bf T}$ is a 2-string tangle with coloring matrix in Eqn.  \ref{eqn_2tangle}.

 \begin{equation}{\bf M_{T}} =
 \left(\begin{array}{c|c}
  A_{(k-2)\times (k-2)}& B_{(k-2)\times 4} \\ \hline
 0_{2\times (k-2)} & \begin{array}{cccc} 
1 & -1 & 1 & -1\\
0 & a & b & c\\
\end{array}
\end{array} \right)
\label{eqn_2tangle}
 \end{equation}

\subsection{Numerator Closure}

\begin{figure}[htbp] \centerline{\psfig{file=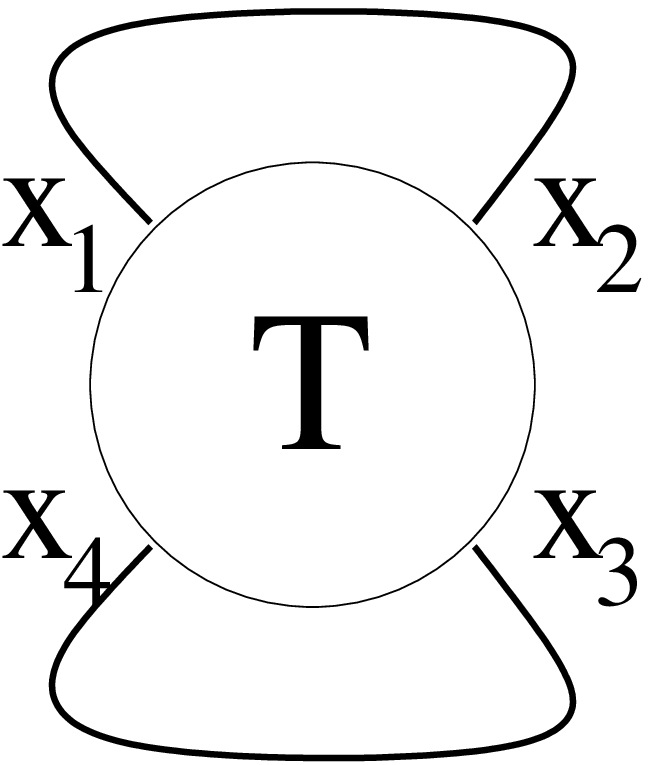,width=.6in}} 
\caption{Numerator Closure, ${\bf N(T)}$}
\label{fig_num}
\end{figure}

The numerator closure of a tangle is formed by connecting the endpoint arcs $x_1$ and $x_2$ as well as $x_3$ and $x_4$ as shown in Fig. \ref{fig_num}.

Thus to obtain a coloring matrix for the knot or link ${\bf N(T)}$, we can add the equations $x_1 - x_2 = 0$ and $x_3 - x_4 = 0$ to a coloring matrix of the tangle ${\bf T}$ (Eqn. \ref{eqn_num2}).  Hence the determinant of the numerator closure of ${\bf T}$,  ${\bf d(N(T))} = |a| {\bf d_U(T)}$.   For example ${\bf d(N({p \over q}))} = |p|$.

\begin{equation}
\left(\begin{array}{c|c}
  A_{(k-2)\times (k-2)}& B_{(k-2)\times 4} \\ \hline
 0_{2\times (k-2)} & \begin{array}{cccc} 
1 & -1 & 1 & -1\\
0 & a & b & c\\
\end{array}
\\ 
 0_{2\times (k-2)} & \begin{array}{cccc} 
1 & -1 & 0 & 0\\
0 & 0 & 1 & -1\\
\end{array}
 \end{array} \right)\sim
\left(\begin{array}{c|c}
  A_{(k-2)\times (k-2)}& B_{(k-2)\times 4} \\ \hline
 0_{4\times (k-2)} & \begin{array}{cccc} 
1 & -1 & 1 & -1\\
0 & a & b & c\\
 0 & 0 & 1 & -1\\
0 & 0 & 0 & 0\\
\end{array}
 \end{array} \right)
\label{eqn_num2}
\end{equation}


\subsection{Denominator Closure}

\begin{figure}[htbp] \centerline{\psfig{file=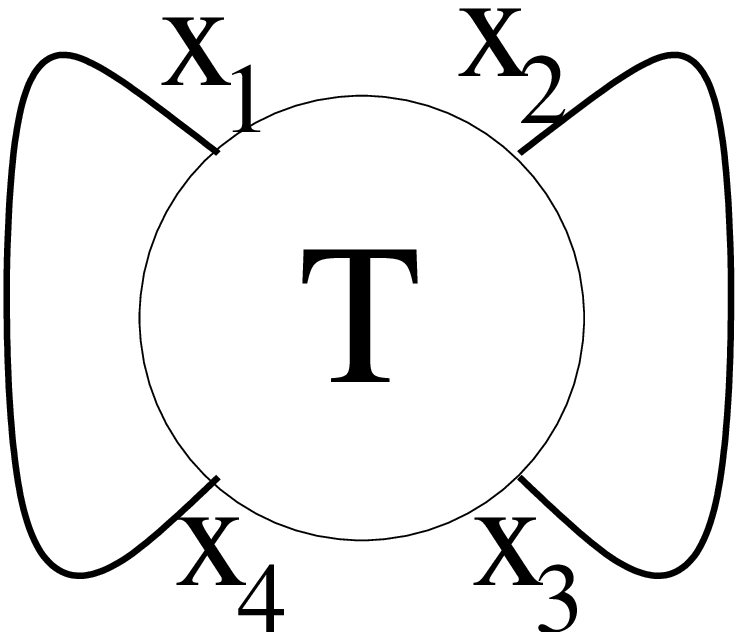,width=.7in}} 
\caption{Denominator Closure, ${\bf D(T)}$}
\label{fig_denom}
\end{figure}

The denominator closure of a tangle is formed by connecting the endpoint arcs $x_1$ and $x_4$ as well as $x_2$ and $x_3$ as shown in Fig. \ref{fig_denom}.

Thus to obtain a coloring matrix for the knot or link ${\bf D(T)}$, we can add the equations $x_1 - x_3 = 0$ and $x_2 - x_3 = 0$ to a coloring matrix of the tangle ${\bf T}$ (Eqn. \ref{eqn_den}).  Hence ${\bf d(D(T))} = |a + b| {\bf d_U(T)}$, 
For example ${\bf d(D({p \over q}))} = |q|$.

\begin{equation}
\left(\begin{array}{c|c}
  A_{(k-2)\times (k-2)}& B_{(k-2)\times 4} \\ \hline
 0_{2\times (k-2)} & \begin{array}{cccc} 
1 & -1 & 1 & -1\\
0 & a & b & c\\
\end{array}
\\ 
 0_{2\times (k-2)} & \begin{array}{cccc} 
1 & 0 & 0 & -1\\
0 & 1 & -1 & 0\\
\end{array}
 \end{array} \right)
\sim
\left(\begin{array}{c|c}
  A_{(k-2)\times (k-2)}& B_{(k-2)\times 4} \\ \hline
 0_{4\times (k-2)} & \begin{array}{cccc} 
1 & -1 & 1 & -1\\
0 & 1 & -1 & 0\\
0 & 0 & a + b & c\\
0 & 0 & 0 & 0\\
\end{array}
 \end{array} \right)
\label{eqn_den}
\end{equation}

\section{Embedding Tangles in Knots/Links}

If a tangle, ${\bf T}$, is a subtangle of a knot/link/tangle, ${\bf K}$, we say that  ${\bf T}$ is embedded in ${\bf K}$.
In the last section, we embedded tangles into knots/links via numerator ${\bf N(T)}$ and denominator ${\bf D(T)}$ closures. 
One can also embed a tangle into a knot/link via much more complicated operations.  Krebes
$\cite{K1}$ proved that if a tangle, ${\bf T}$, is a subtangle of a knot or link, ${\bf K}$, then $gcd ({\bf d(N(T), D(T)})$ divides ${\bf d(K)}$.  A short proof of this result is also given in $\cite{R1, SWK}$.  We will also provide a short proof of this result.
A similar technique to that presented here was used in \cite{P1} to prove several related results.

Recall that if ${\bf M_l(T)}$ is row equivalent to 
$\left(\begin{array}{cccc} 
1 & -1 & 1 & -1\\
0 & a & b & c\\
\end{array}\right)$, then ${\bf d(N(T))} = |a|{\bf d_U(T)}$ and ${\bf d(D(T))} = |b+a|{\bf d_U(T)}$.  Hence $gcd {\bf(d(N(T), D(T))} = gcd(|a|{\bf d_U(T)}, |a+b|{\bf d_U(T)}) = {\bf d_U(T)}gcd(a, b)$. 

Let $g = gcd(a, b)$.  Since we started out with a matrix where the sum of the entries in a row is 0, $a + b + c = 0$.  Thus $g$ also divides $c$.  Let $a' = {a \over g}$,  $b' = {b \over g}$,  $c' = {c \over g}$.  Suppose ${\bf T}$ is embedded in a knot ${\bf K}$.  Then the matrix $M_1$ in equation \ref{krebes} is a coloring matrix of ${\bf K}$ where the upper left $k \times (k+2)$ submatrix is a coloring matrix of ${\bf T}$. The matrix $M_2$ in equation \ref{krebes} is obtained form the matrix $M_1$ by dividing a row by $g$.  Hence $det(M_1) = g det (M_2)$.  Since ${\bf d_U(T)} = det(A)$ and $det(A)$ divides $det (M_2)$,  ${\bf d_U(T)}gcd(a, b)$ divides $det(M_1)$.  Thus $gcd {\bf(d(N(T), D(T))} $ divides ${\bf d(K)}$.

\begin{equation}
M_1 = \left(
\begin{array}{c|c|c}
  	\begin{array}{c}
   		A_{(k-2)\times (k-2)} \\ \hline
   		0_{2\times (k-2)} \cr
		~  \cr
  	\end{array}
&
	\begin{array}{c}
   		B_{(k-2)\times 4} \\ \hline
   		\begin{array}{cccc} 
   			1 & -1 & 1 & -1\\
   			0 & a & b & c\\
   		\end{array}
	\end{array}
& ~~  0_{k\times (c-2)}  \cr  \hline
0_{c\times (k-2)}  & D_{c \times 4}  & ~~E _{c \times (c-2)}
\end{array}
\right)\hfill  
~~~~  M_2 = \left(
\begin{array}{c|c|c}
  	\begin{array}{c}
   		A_{(k-2)\times (k-2)} \\ \hline
   		0_{2\times (k-2)} \cr
		~  \cr
  	\end{array}
&
	\begin{array}{c}
   		B_{(k-2)\times 4} \\ \hline
   		\begin{array}{cccc} 
   			1 & -1 & 1 & -1\\
   			0 & a' & b' & c'\\
   		\end{array}
	\end{array}
& ~~  0_{k\times (c-2)}  \cr  \hline
0_{c\times (k-2)}  & D_{c \times 4}  & ~~E _{c \times (c-2)}
\end{array}
\right)
\label{krebes}
\end{equation}

\section{How good of a tangle invariant is
colorability?}

\begin{figure}[htbp] 
{ \vskip 13pt}
\centerline{\psfig{file=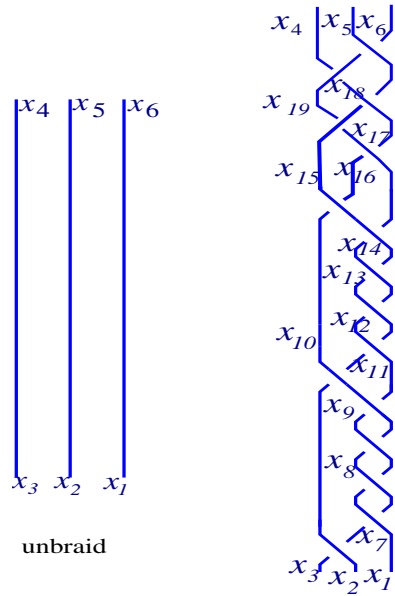,width=2.0in,height=3.1in}}
{ \vskip 13pt} 
\caption{Coloring cannot distinguish the unbraid from other
braids.} \label{fig:unb}
\end{figure}

Recall the ${\bf d_U(B)}$ = 1 for all braids ${\bf B}$. For the unbraid, ${\bf U}$, shown on the left in Fig. $\ref{fig:unb}$, ${\bf M_{l}(U)}$ is given in Eqn. \ref{eqn:ub}.  This invariant is the same for the braid shown on the right-side of Fig. $\ref{fig:unb}$.  Thus the coloring invariants are the same for these two braids.  Hence,
coloring cannot distinguish the unbraid from all other braids.

\begin{equation}{\bf M_{l}(U)} =
 \left(\begin{array}{cccccc}
 1 & 0 & 0& 0 & 0& -1\\
0 & 1 & 0 & 0 & -1 & 0\\
0 & 0 & 1 & -1 &0 & 0\\
\end{array} \right)
\label{eqn:ub}
\end{equation}

\subsection{ Coloring Can Distinguish
Between A Tangle And Its Mirror Image}
\label{subs:mirr}

Coloring mod $m$ cannot distinguish between a knot and
its mirror image. If a knot is $m$-colorable, then so is its mirror image.  For example  the trefoil knot, ${\bf N({3 \over 1})}$, and its mirror image, ${\bf N({3 \over -1})}$, are both
3-colorable.  However, coloring can distinguish between the rational tangle ${\bf {3 \over 1} }$ and its mirror image ${\bf {3 \over -1} }$:
${\bf M_l({3 \over 1})} = 
\left( \begin{matrix}
1 ~~& -1 ~~& 1      ~~& -1   \cr
0 ~~& 3  ~~& -4 ~~& 1
\end{matrix} \right)$ while ${\bf M_l({3 \over -1})} = 
\left( \begin{matrix}
1 ~~& -1 ~~& 1      ~~& -1   \cr
0 ~~& 3  ~~& -2 ~~& 1
\end{matrix} \right)$.

\baselineskip 4.8mm
   \bibliography{../mybib}
    \bibliographystyle{plain}

\end{document}